\documentclass{article}
\usepackage{amsmath,amssymb,latexsym,amsthm,amscd}
\font\goth=eusm10

\newcommand\bc{\mathbf C}

\newcommand\E{\mathcal E}

\newcommand\Ii{\hbox{\goth I}}

\newcommand\bz{\mathbf{Z}}

\newcommand\bn{\mathbf{N}}

\newcommand\Oc{\hbox{\goth O}}

\newtheorem{theorem}{Theorem}

\newtheorem{lemma}{Lemma}

\numberwithin{proposition}{section}
\numberwithin{definition}{section}
\numberwithin{corollary}{section}
\numberwithin{remark}{section}
\numberwithin{lemma}{section}
\numberwithin{equation}{section}

\numberwithin{theorem}{section}

\numberwithin{question}{section}
\numberwithin{case}{section}
\numberwithin{example}{section}
\numberwithin{conjecture}{section}

\begin{document}
\title{On a classical correspondence between \\ K3 surfaces III} 
\author{C.G.Madonna \footnote{Supported by I3P contract} \ and 
Viacheslav V.Nikulin \footnote{Supported by EPSRC grant  
EP/D061997/1}}
\maketitle
\begin{abstract} Let $X$ be a K3 surface, and $H$ its primitive polarization 
of the degree $H^2=8$. The moduli space of sheaves over $X$ with the 
isotropic Mukai vector $(2,H,2)$ is again a K3 surface, $Y$. In 
\cite{Mad-Nik1} we gave necessary and sufficient conditions in terms of 
Picard lattice of $X$ when $Y$ is isomorphic to $X$. 
The proof of sufficient condition in \cite{Mad-Nik1}, 
when $Y$ is isomorphic to $X$, used Global Torelli Theorem for K3 surfaces, 
and it was not effective. 

Here we give an effective variant of these results: its sufficient 
part  gives an explicit isomorphism between $Y$ and $X$.   

We hope that our similar results in \cite{Mad-Nik2}, \cite{Nik1}, 
\cite{Nik2} for arbitrary primitive isotropic Mukai vector on a K3 surface
also can be  made effective.  
\end{abstract}

\section{Introduction} \label{S:intro}

In \cite{Mad-Nik1} we had obtained the following result.

\begin{theorem} \label{maintheorem1}
Let $X$ be a K3 surface over $\bc$ with Picard lattice $N(X)$, 
and $H\in N(X)$ is primitive, nef with $H^2=8$. 
Let $Y$ be the moduli space of sheaves 
on $X$ with the isotropic Mukai vector $v=(2,H,2)$. 

Then, $Y\cong X$ if there exists $h_1\in N(X)$ such that the primitive 
sublattice $[H,h_1]_{\text{pr}}$ in $N(X)$ generated by $H$ and $h_1$ has 
an odd determinant (equivalently, $H\cdot [H,h_1]_{\text{pr}}=\bz$)
and
$$
h_1^2=\pm 4\ \ \ and \ \ h_1\cdot H\equiv 0\mod 2.
$$

These conditions are necessary for $Y\cong X$ if the Picard number 
$\rho(X)=\text{rk}\ N(X)$ $\le 2$, and $X$ is a general K3 surface with its Picard lattice 
(i.e. the automorphism group of the transcendental periods 
$(T(X), H^{2,0}(X))$ is $\pm 1$). 
\end{theorem}

The proof of Theorem \ref{maintheorem1} given in \cite{Mad-Nik1} used 
Global Torelli Theorem for K3 surfaces \cite{PShShaf}, and it 
was not effective; under conditions of Theorem \ref{maintheorem1}, 
we had only proved existence of the isomorphism $Y\cong X$.  

The purpose of this paper is to proof the following effective 
variant of Theorem \ref{maintheorem1}.  

\begin{theorem} 
\label{maintheoremef1}
Let $X$ be a K3 surface over $\bc$ and $H\in N(X)$ 
is primitive, nef with $H^2=8$. Let $Y$ be the moduli space of sheaves 
on $X$ with the isotropic Mukai vector $v=(2,H,2)$. 

Then $Y\cong X$ with an explicit geometric isomorphism given by 
\eqref{isomcase1} and \eqref{isomcase2} in the proof 
below if there exists $h_1\in N(X)$ such that the primitive 
sublattice $[H,h_1]_{\text{pr}}$ in $N(X)$ generated by $H$ and $h_1$ has 
an odd determinant (equivalently, $H\cdot [H,h_1]_{\text{pr}}=\bz$), 
and
\begin{equation}
h_1^2=\pm 4,\ \ h_1\cdot H\equiv 0\mod 2,
\label{condmain-1}
\end{equation}
and 
\begin{equation}
h^0\Oc_X(h_1)=h^0\Oc_X(-h_1)=0\ \ if\ \ h_1^2=-4.
\label{condmain0}
\end{equation}

These conditions are necessary for $Y\cong X$ if the Picard number 
$\rho(X)\le 2$, and $X$ is a general K3 surface with its Picard lattice 
(i.e. the automorphism group of the transcendental periods 
$(T(X), H^{2,0}(X))$ is $\pm 1$). 
\end{theorem}

Here results by Tyurin \cite{Tyurin1} and by Ballico-Chiantini 
\cite{BC} are very useful. 

See Remark 3.1 about difference between 
conditions of Theorems \ref{maintheorem1} and \ref{maintheoremef1}. 

\par\medskip

In our papers \cite{Mad-Nik2}, \cite{Nik1}, \cite{Nik2}, Theorem 
\ref{maintheorem1} was generalized to arbitrary primitive Mukai 
vector. We hope that similar considerations as in this paper will also  
permit to make these results effective. 
\par\medskip

The first author thanks Profs L. Chiantini, A. Verra, and 
Dr. A. Rapagnetta for useful discussions.

\section{Reminding of the Main Result of \cite{Mad-Nik1}}
\label{section2}

We denote by $X$ an algebraic K3 surface 
over the field $\bc$ of complex numbers. I.e. $X$ is a non-singular 
projective algebraic surface over $\bc$ with the trivial canonical class 
$K_X = 0$ and the vanishing  irregularity $q(X)=0$. 

We denote by $N(X)$ 
the Picard lattice (i.e. the lattice of 2-dimensional algebraic 
cycles) of $X$. By $\rho (X)=\text{rk}\ N(X)$ we denote the Picard number of 
$X$. By 
\begin{equation}
T(X)=N(X)_{H^2(X,\bz)}^\perp 
\end{equation}
we denote the 
transcendental lattice of $X$.

\par\medskip

For a Mukai vector 
$v=(r,c_1,s)$ where $r,s\in \bz$ and $c_1\in N(X)$, 
we denote by $Y=M_X(r,c_1,s)$
the moduli space of stable (with respect to some ample 
$H^\prime\in N(X)$) rank $r$ sheaves on $X$ with first Chern classes $c_1$, and Euler characteristic $r+s$.

By results of Mukai \cite{Muk1}, \cite{Mukai2},
under suitable conditions on the Chern classes, the moduli space $Y$ is 
always deformations equivalent
to a Hilbert scheme of 0-dimensional cycles on $X$ (of same dimension).
\par\medskip

In  \cite{Mad-Nik1} we had considered the case of the isotropic 
Mukai vector $v=(2,H,2)$ with $H^2=8$ and $H$ nef and primitive, and
we had looked for conditions on the Picard lattice $N(X)$ which imply 
that $Y\cong X$. One of our main results in \cite{Mad-Nik1} was 
the following Theorem.

\begin{theorem} Let $X$ be a K3 surface over $\bc$ and $H\in N(X)$ 
is primitive, nef with $H^2=8$. Let $Y$ be the moduli space of sheaves 
on $X$ with the isotropic Mukai vector $v=(2,H,2)$. 

Then, $Y\cong X$ if there exists $h_1\in N(X)$ such that the primitive 
sublattice $[H,h_1]_{\text{pr}}$ in $N(X)$ generated by $H$ and $h_1$ has 
an odd determinant (equivalently, $H\cdot [H,h_1]_{\text{pr}}=\bz$)
and
$$
h_1^2=\pm 4\ \ \ and \ \ h_1\cdot H\equiv 0\mod 2.
$$

These conditions are necessary for $Y\cong X$ if the Picard number 
$\rho(X)\le 2$, and $X$ is a general K3 surface with its Picard lattice 
(i.e. the automorphism group of the transcendental periods 
$(T(X), H^{2,0}(X))$ is $\pm 1$). 
\label{maintheorem}
\end{theorem}

The proof of Theorem \ref{maintheorem} in \cite{Mad-Nik1} used 
Global Torelli Theorem for K3 surfaces \cite{PShShaf}, and it 
was not effective. 

The purpose of this paper is to give an effective variant of Theorem 
\ref{maintheorem} which does not use Global Torelli Theorem 
for K3 surfaces. This variant gives an explicit isomorphism $Y\cong X$. In 
\cite{Mad-Nik1} we had only proved existence of such isomorphism. 

\section{An effective (without use of Global Torelli Theorem for K3) 
variant of Theorem \ref{maintheorem}}\label{section3}

Here we prove the following effective variant of Theorem \ref{maintheorem}.

\begin{theorem} Let $X$ be a K3 surface over $\bc$ and $H\in N(X)$ 
is primitive, nef with $H^2=8$. Let $Y$ be the moduli space of sheaves 
on $X$ with the isotropic Mukai vector $v=(2,H,2)$. 

Then $Y\cong X$ with an explicit geometric isomorphism given by 
\eqref{isomcase1} and \eqref{isomcase2} in the proof 
below if there exists $h_1\in N(X)$ such that the primitive 
sublattice $[H,h_1]_{\text{pr}}$ in $N(X)$ generated by $H$ and $h_1$ has 
an odd determinant (equivalently, $H\cdot [H,h_1]_{\text{pr}}=\bz$), 
and
\begin{equation}
h_1^2=\pm 4,\ \ h_1\cdot H\equiv 0\mod 2,
\label{condmain-1}
\end{equation}
and 
\begin{equation}
h^0\Oc_X(h_1)=h^0\Oc_X(-h_1)=0\ \ if\ \ h_1^2=-4.
\label{condmain0}
\end{equation}

These conditions are necessary for $Y\cong X$ if the Picard number 
$\rho(X)\le 2$, and $X$ is a general K3 surface with its Picard lattice 
(i.e. the automorphism group of the transcendental periods 
$(T(X), H^{2,0}(X))$ is $\pm 1$). 
\label{maintheoremef}
\end{theorem}

\begin{proof}
The `necessary' part of the proof is the same as in \cite{Mad-Nik1}.  
Let us assume that $\rho(X)\le 2$, $X$ is general 
(i. e. the automorphism group of the transcendental 
periods $(T(X),H^{2,0}(X))$ is $\pm 1$), and $Y\cong X$. 
Since $Y\cong X$, periods of $Y$ and $X$ must be isomorphic.  
We have shown in \cite{Mad-Nik1} that periods of $Y$ and $X$ 
are isomorphic if and only if $H\cdot N(X)=\bz$ (Mukai 
condition), and there exists $h_1\in N(X)$ 
such that  
$$
h_1^2=\pm 4,\ \  H\cdot h_1\equiv 0\mod 2.
$$
Thus, we obtain exactly the conditions of of Theorem \ref{maintheorem}, or
conditions of Theorem \ref{maintheoremef} except \eqref{condmain0}. 
The determinant of 
the Gram 
matrix of $H$ and $h_1$ is equal to 
\begin{equation}
H^2h_1^2-(H\cdot h_1)^2=\pm 32-(H\cdot h_1)^2\not =0.
\label{eqrank2}
\end{equation}
Thus, $\rho (X)=2$, and $N(X)=[H,h_1]_{\text{pr}}$ is a 2-dimensional lattice 
($\rho(X)=1$ never happens if $Y\cong X$). 
In (\cite{Mad-Nik1}, 
Proposition 3.2.1 and Theorems 3.2.2, 3.2.3) we have shown that if 
$N(X)$ has only elements $h_1$ satisfying \eqref{condmain-1}  
with $h_1^2=-4$, then $N(X)$ has no elements $\delta$ with $\delta^2=-2$. 
Since any irreducible curve $C$ on a K3 surface has $C^2\ge -2$ (it is 
well-known and obvious) and $N(X)$ is an even lattice, it 
then follows that any effective element of $N(X)$ has a non-negative 
square. Then \eqref{condmain0} is automatically valid. 
 
Now let us consider the `sufficient' part of the proof of Theorem  
\ref{maintheoremef} which used 
Global Torelli Theorem for K3 surfaces \cite{PShShaf} and was not 
effective in \cite{Mad-Nik1}. 

We have simple 

\begin{lemma} Let $X$ be a K3 surface and $H\in N(X)$ a primitive element 
with $H^2=8$.

Then existence of $h_1\in N(X)$ satisfying conditions of Theorem 
\ref{maintheorem}, i. e. 
\begin{equation}
\left(\ h_1^2=\pm 4,\ H\cdot h_1\equiv 0\mod 2,\ H\cdot[H,h_1]_{\text{pr}}=\bz 
\ \right)
\label{condmain1}
\end{equation}
is equivalent to 
\begin{equation}
\exists D\in N(X)\ such\ that\ Mukai\ vector\ v_1=(2,H+2D,\pm 1)\ is\ isotropic
\label{condmain2}
\end{equation}
i.e. $(H+2D)^2=\pm 4$. 
\par\medskip
The relation between \eqref{condmain1} and \eqref{condmain2} is just 
\begin{equation}
h_1=H+2D\,.
\label{condmain3}
\end{equation}

\label{mainlemma}
\end{lemma}

\begin{proof} Assume \eqref{condmain1} is valid. The determinant of 
Gram matrix of $H$ and $h_1$ is equal to $\pm 32-(H\cdot h_1)^2\not=0$. 
It follows that $[H,h_1]_{\text{pr}}$ is a 2-dimensional sublattice in 
$N(X)$. 
 
Since $H$ is primitive in $N(X)$, then $H\notin 2[H,h_1]_{\text{pr}}$. 
Since $(h_1/2)^2=\pm 1$ and $N(X)$ is even lattice, then 
$h_1\notin 2[H,h_1]_{\text{pr}}$. If 
$H-h_1\notin 2[H,h_1]_{\text{pr}}$, it then follows that 
$H+2[H,h_1]_{\text{pr}}$, $h_1+2[H,h_1]_{\text{pr}}$ give a 
basis of 
$[H,h_1]_{\text{pr}}\mod 2=[H,h_1]_{\text{pr}}/2[H,h_1]_{\text{pr}}$. 
Then $H\cdot [H,h_1]_{\text{pr}}\equiv \{H^2,\ H\cdot h_1\}\equiv 0\mod 2$ 
which contradicts $H\cdot [H,h_1]_{\text{pr}}=\bz$. Thus, 
$h_1=H+2D$ where $D\in [H,h_1]_{\text{pr}}\subset N(X)$. It follows 
the condition \eqref{condmain2}. 

Now assume \eqref{condmain2} is valid. We put $h_1=H+2D$. Then 
$h_1^2=(H+2D)^2=\pm 4$, $h_1\cdot H=H^2+2(H\cdot D)\equiv 0\mod 2$. 
We have 
$$
(H+2D)^2=8+4(H\cdot D)+4D^2=\pm 4
$$
where $D^2\equiv 0\mod 2$ since $N(X)$ is even. It follows  
$H\cdot D\equiv 1\mod 2$. Since $H\cdot H=8$, it follows 
$H\cdot [H,h_1]_{\text{pr}}=\bz$. We obtain the condition 
\eqref{condmain1}. 

This finishes the proof of Lemma \ref{mainlemma}. 
\end{proof}

\medskip

To give an effective proof of Theorem \ref{maintheoremef}, we now should 
consider two cases. 

\medskip 

{\it The case $h_1^2=4$ of Theorem \ref{maintheoremef}.}
By Lemma \ref{mainlemma}, this is equivalent to the  
\begin{equation}
\text{Mukai vector\ } 
v_1=(2,h_1=H+2D,1)\ \text{is isotropic for some} \ D\in N(X).
\label{condmain4}
\end{equation}
Then $Y=M_X(v) \cong M_X(v_1)$ under tensorization by $\Oc_X(D)$.
By general results (see e.g. \cite{Tyurin1}, Chapter II, Section 4)  
$M_X(v_1) \cong M_X(w_1) \cong X$ where
$w_1=(1,h_1,2)$. This gives an explicit isomorphism 
\begin{equation}
Y=M_X(2,H,2)\cong M_X(2,H+2D,1)\cong M_X(1,H+2D,2)\cong X\,. 
\label{isomcase1}
\end{equation}

\medskip

{\it The case $h_1^2=-4$ of Theorem \ref{maintheoremef}.}
By Lemma \ref{mainlemma}, this is equivalent to the  
\begin{equation}
\text{Mukai vector\ }v_1=(2,h_1=H+2D,-1)\ \text{is isotropic for some} 
\ D\in N(X),
\label{condmain5}
\end{equation}
and
\begin{equation}
h^0\Oc_X(H+2D)=h^0\Oc_X(-H-2D)=0.
\label{condmain6}
\end{equation}

Let us take $D\in N(X)$ satisfying these conditions. 
Changing $h_1$ by $-h_1$ is equivalent to changing $D$ by $-H-D$. 
Replacing $h_1$ by $-h_1$ if necessary, we can assume 
that $H\cdot h_1=H\cdot (H+2D)=8+2H\cdot D\ge 0$. Equivalently, 
$H\cdot D\ge -4$. From $(H+2D)^2=8+4H\cdot D+4D^2=-4$ and $D^2\equiv 0\mod 2$, 
it follows that $H\cdot D$ is always odd (and $h_1\cdot D$ as well). 
Thus, we can even assume more:  
\begin{equation}
 H\cdot D>-4\,.
\label{cond6}
\end{equation}

\par\medskip

From $h^0 \Oc_X(h_1)=0$ and $h^2\Oc_X(h_1)=h^0 \Oc_X(-h_1) =0$ and 
Riemann-Roch Theorem for K3, we obtain that $\chi \Oc_X(h_1)=0$ and 
$h^1\Oc_X(h_1)=0$. 

\par\medskip

Let $p\in X$ be a point and $\Ii_p$ its sheaf of ideals. Since 
$\Ii_p\subset \Oc_X$ and $h^0\Oc_X(h_1)=0$, 
then $h^0 \Ii_p(h_1)=0$. By
the exact sequence of $\Ii_p \subset \Oc_X$, we also obtain
$h^1 \Ii_p(h_1)=h^1\Ii_p(H+2D)=h^0(\Oc_p(h_1))=1$. 
Then $(H^1\Ii_p(H+2D))^\ast\cong Ext^1(\Ii_p(H+D),\Oc_X(-D))$ is 
one-dimensional, and we can 
construct a rank 2 bundle $\E$ given by
the non-trivial extension
\begin{equation}
0 \to \Oc_X(-D) \to \E \to \Ii_p(H+D) \to 0, 
\label{exact2}
\end{equation}
equivalently a rank 2 bundle $\E(D)$ 
\begin{equation}
0 \to \Oc_X \to \E(D) \to \Ii_p(H+2D) \to 0,
\label{exact1}
\end{equation}
and $\E$ is a rank 2 bundle with $c_1=H$ and $c_2=4$. The bundle $\E$ 
is semistable since $\E(D)$ is so.
If $L=\Oc_X(L) \subset \E(D)$ is such that $L\cdot H>(H \cdot h_1)/2>0$ then $L$ is not contained
in the image of the map $\Oc_X \to \E(D)$.
Hence the image of the inclusion of $L \subset \E(D)$ under the projection
$\E \to \Ii_p(H+2D)$ gives a 
non zero map $L \to \Ii_p \otimes \Oc_X(H+2D)$ and $h_1=L+L'$ with 
$L'$ effective. 
Then $h_1-L=L'$ is effective and we have
$h^0 \Oc_X(h_1-L) \leq h^0 \Oc_X(h_1)=0$ which is absurd.
Indeed the last vanishing follows by the exact sequence 
\begin{equation}
0 \to \Oc_X(h_1-L) \to \Oc_X(h_1) \to \Oc_L(h_1) \to 0
\end{equation}
Thus $\E \in M_X(2,H,2)$. 

Since $h^0\Ii_p(H+2D)=0$, by  \eqref{exact1}  
we obtain that $h^0\E(D)=1$. Thus, \eqref{exact2} is defined by 
a unique (up to proportionality) non-zero section of $\E(D)$, 
and $p$ is the zero locus of this section. 
Thus the constructed using \eqref{exact1} and \eqref{exact2} map 
\begin{equation}
X\to M_X(2,H,2)=Y
\label{isomcase2}
\end{equation} 
has the degree one, and it defines an explicit isomorphism $Y \cong X$. 

\medskip 

This finishes the proof of Theorem \ref{maintheoremef}. 
\end{proof}

\noindent
{\bf Remark 3.1.}  
Difference between Theorems \ref{maintheorem} and 
\ref{maintheoremef} is in condition \eqref{condmain0} which means 
that both elements $\pm h_1$ with $h_1^2=-4$ should be also not effective. 

If there exists $h_1^\prime\in N(X)$ with $(h_1^\prime)^2=-4$, then 
$N(X)$ has plenty of elements $h_1$ with $h_1^2=-4$ such that both 
$\pm h_1$ are not effective. 

Really, the nef cone $NEF(X)$ of $X$ is a fundamental 
chamber for the group $W^{(-2)}(X)$ generated by reflections in all elements 
$\delta\in N(X)$ with $\delta^2=-2$. It follows that there exists 
$w\in W^{(-2)}(X)$ such that $h_1=w(h_1^\prime)$ divides $NEF(X)$ in two 
open parts: there exist two nef elements $H_1, H_2 \in NEF(X)$ such that 
$H_1\cdot h_1<0$ and $H_2\cdot h_1>0$ (such elements $\pm h_1$ are 
called {\it not pseudo-effective}). It follows that both elements 
$\pm h_1$ are not effective. 

Thus, to satisfy conditions of Theorem \ref{maintheoremef} for  
the case $h_1^2=-4$,  
one should look first for not pseudo-effective $\pm h_1$ satisfying 
$h_1^2=-4$, $h_1\cdot H\equiv 0\mod 2$ and $h_1\cdot [H,h_1]_{\text{pr}}=\bz$. 
All not pseudo-effective elements $\pm f\in N(X)$ with negative square $f^2<0$ 
satisfy the geometric condition \eqref{condmain0}, and there are plenty 
of them. 

\medskip 

In our papers \cite{Mad-Nik2}, \cite{Nik1}, \cite{Nik2}, Theorem 
\ref{maintheorem} was generalized to arbitrary primitive Mukai 
vector $v=(r,H,s)$ where $r,s\in \bn$ and $H^2=2rs$. 
We hope that similar considerations as here also permit 
to make these results effective. We hope to consider that in 
further publications.

\

\

C.G.Madonna \par
Math. Dept., 
CSIC, C/ Serrano 121, 
28006 Madrid,
SPAIN

carlo@madonna.rm.it \ \
cgm@imaff.cfmac.csic.es

\

\

V.V.Nikulin \par
Deptm. of Pure Mathem. The University of Liverpool, Liverpool\par
L69 3BX, UK;
\vskip1pt
Steklov Mathematical Institute,\par
ul. Gubkina 8, Moscow 117966, GSP-1, Russia

vnikulin@liv.ac.uk \ \
vvnikulin@list.ru

\end{document}